\documentclass[10pt]{article}
\usepackage{latexsym,amsfonts,amssymb}
\begin{document}

\title{ {\bf Stochastic quasi-geostrophic equation
\thanks{Research supported by 973 project, NSFC, key Lab of CAS, the DFG through IRTG 1132 and CRC 701}\\} }
\author{ {\bf Michael R\"{o}ckner}$^{\mbox{a}}$, {\bf Rongchan Zhu}$^{\mbox{a,b},}$, {\bf Xiangchan Zhu}$^{\mbox{a,c},}$
\thanks{E-mail address: roeckner@math.uni-bielefeld.de(M. R\"{o}ckner),
zhurongchan@126.com(R. C. Zhu), zhuxiangchan@126.com(X. C. Zhu)}\\ \\
$^{\mbox{a}}$Department of Mathematics, University of Bielefeld, D-33615 Bielefeld, Germany,\\
$^{\mbox{b}}$Institute of  Applied Mathematics, Academy of
Mathematics and
Systems Science, \\
Chinese Academy of Sciences, Beijing 100190, China,\\
$^{\mbox{c}}$School of Mathematical Sciences, Peking University, Beijing 100871, China}
\date{}
\maketitle

 \noindent\mbox{}\hrulefill\mbox{}

\noindent {\bf Abstract}

In this note we study the 2d stochastic quasi-geostrophic equation in $\mathbb{T}^2$ for general parameter $\alpha\in (0,1)$ and multiplicative noise. We prove the existence of martingale solutions and pathwise uniqueness under some condition in the general case , i.e. for all $\alpha\in (0,1)$ .
In the subcritical case $\alpha>1/2$, we prove  existence and uniqueness of (probabilistically) strong solutions and construct a Markov family of solutions. In particular, it is uniquely ergodic for $\alpha>\frac{2}{3}$ provided the noise is non-degenerate. In this case, the convergence to the (unique) invariant measure is exponentially fast. In the general case, we prove the existence of Markov selections.

\noindent\mbox{}\hrulefill\mbox{}

\renewcommand{\theequation}{1.\arabic{equation}}

\noindent{\bf 1. Introduction and notation}------Consider the following two dimensional (2D) stochastic quasi-geostrophic equation in the periodic domain $\mathbb{T}^2=\mathbb{R}^2/(2\pi \mathbb{Z})^2$:
$$\frac{\partial \theta(t,x)}{\partial t}=-u(t,x)\cdot \nabla \theta(t,x)-\kappa (-\triangle)^\alpha \theta(t,x)+G(\theta,\xi)(t,x),\eqno(1.1)$$
with initial condition $\theta(0,x)=\theta_0(x),$
where $\theta(t,x)$ is a real-valued function of $x$ and $t$, $0<\alpha<1, \kappa>0$ are real numbers. $u$ is determined by $\theta$ through a stream function $\psi$ via the following relations:
$$u=(u_1,u_2)=(-R_2\theta,R_1\theta).\eqno(1.2)$$
Here $R_j$ is the $j$-th periodic Riesz transform and $\xi(t,x)$ is a Gaussian random field, white noise in time, subject to the restrictions imposed below. The case $\alpha=\frac{1}{2}$ is called the critical case, the case $\alpha>\frac{1}{2}$ sub-critical and the case $\alpha<\frac{1}{2}$ super-critical. The existence of weak solutions in the deterministic case has been obtained in \cite{Re95}. In the following, we will restrict ourselves to flows which have zero average on the torus, i.e.$\int_{\mathbb{T}^2}\theta dx=0.$
Set $H=L^2(\mathbb{T}^2)$ and let $|\cdot|$ and $\langle .,.\rangle$ denote the norm and inner product in $H$, respectively. We recall that on $\mathbb{T}^2$, $\sin (k\cdot x),\cos(k\cdot x)$ form an eigenbasis of $-\triangle$. Here $k\in \mathbb{Z}^2\backslash \{0\}, x\in \mathbb{T}^2$ and the corresponding eigenvalues are $|k|^2$. Define $\|f\|_{H^s}^2:=\sum_k |k|^{2s}\langle f,e_k\rangle^2$ and let $H^s$ denote the Sobolev space of all $f$ for which $\|f\|_{H^s}$ is finite. Set $\Lambda=(-\triangle)^{1/2}$.
Define the linear operator $A: D(A)\subset H\rightarrow H$ as $Au=\kappa (-\triangle)^\alpha u.$ The operator $A$ is positive definite and selfadjoint. Denote the eigenvalues of $A$ by $0<\lambda_1\leq\lambda_2\leq\cdots$ , and by $e_1,e_2,$... the corresponding complete orthonormal system in $H$ of eigenvectors of $A$. We also denote $\|u\|=|A^{1/2}u|$, then $\|\theta\|^2\geq\lambda_1|\theta|^2$.

\noindent{\bf 2. Existence and uniqueness of solutions }-----By the above definitions  Eqs (1.1)-(1.2) turn into the abstract stochastic evolution equation
$$\left\{\begin{array}{ll}d\theta(t)+A\theta(t)dt+u(t)\cdot \nabla\theta(t)dt= G(\theta(t))dW(t),&\ \ \ \ \textrm{ }\\\theta(0)=x,&\ \ \ \ \textrm{ } \end{array}\right.\eqno(2.1)$$
where $u$ satisfies (1.2) and $W(t)$ is a cylindrical Wiener process in a separable Hibert space $K$ defined on a probability space $(\Omega,\mathcal{F},P)$. Here $G$ is a mapping from $H^\alpha$ to $L^2(K,H)$.

\noindent{\bf Definition 2.1} We say that there exists a martingale solution to (2.1) if there exists a stochastic basis $(\Omega,\mathcal{F},\{\mathcal{F}_t\}_{t\in [0,T]},P)$, a cylindrical Wiener process $W$ on the space $K$ and a progressively measurable process $\theta:[0,T]\times \Omega\rightarrow H$, such that for $P$-a.e $\omega\in \Omega$

$$\theta(\cdot,\omega)\in L^\infty(0,T;H)\cap L^2(0,T;H^\alpha)\cap C([0,T];H_w)\eqno(2.2)$$and $P$-a.s.
$$\langle \theta(t),\psi\rangle+\int_0^t\langle A^{1/2}\theta(s),A^{1/2}\psi\rangle ds-\int_0^t \langle u(s)\cdot \nabla \psi,\theta(s)\rangle ds=\langle x,\psi\rangle+\langle \int_0^tG(\theta(s))dW(s),\psi\rangle, \eqno(2.3)$$  for all $t\in [0,T]$ and all $\psi\in C^1(\mathbb{T}^2)$.
Here $C([0,T];H_w)$ denotes the space of $H$-valued weakly continuous functions on $[0, T]$.

\noindent{\bf Remark 2.2} Note that for regular functions $\theta$ and $v$, we have
$\langle u(s)\cdot \nabla (\theta(s)+\psi),\theta(s)+\psi\rangle=0,$ so
$\langle u(s)\cdot \nabla \theta(s),\psi\rangle=-\langle u(s)\cdot \nabla \psi,\theta(s)\rangle.$ Thus the integral equation (2.3) corresponds to equation (2.1).

\noindent{\bf Definition 2.3} We say that there exists a (probabilistically strong) solution to (2.1) over the time interval $[0,T]$ if for every probability space $(\Omega,\mathcal{F},\{\mathcal{F}_t\}_{t\in [0,T]},P)$ with an $\mathcal{F}_t$-Wiener process $W$, there exists a a progressively measurable process $\theta:[0,T]\times \Omega\rightarrow H$ such that (2.2) and (2.3) hold.

\noindent{\bf 2.1. The general case}-----Consider the following condition:

\noindent $(G.1) \textrm{ }G:H\rightarrow L_2(K,H)\textrm{ is continuous and }|G(\theta)|^2_{L_2(K,H)}\leq \lambda_0|\theta|^2+\rho,\theta\in H$ for some positive real numbers $\lambda_0$ and $\rho$.\\
By the compactness method based on fractional Sobolev spaces in \cite{FG95}, we obtain the existence of martingale solutions.\\
\noindent{\bf Theorem 2.1.1} Under condition  (G.1), there is a martingale solution $(\Omega,\mathcal{F},\{\mathcal{F}_t\},P,W, \theta)$ to (2.1). Moreover, if $x\in H^1$ and $G$ satisfies $$\int(\sum_j|G(\theta)(e_j)|^2)^{p/2}dx\leq C(1+\int |\theta|^pdx), \forall t>0\eqno(2.4)$$ with $2<p<\infty$, then
we have $\sup_{t\in [0,T]}\|\theta(t)\|_{L^p}<\infty,P-a.s.$

\noindent{\bf Theorem 2.1.2} If $G\in L^2(K,H)$ does not depend on $\theta$, then there is a martingale solution $(\Omega,\mathcal{F},\{\mathcal{F}_t\},P,W, \theta)$ to (2.1). Moreover, if $x\in L^p(\mathbb{T}^2)$, then
 $\sup_{t\in [0,T]}\|\theta(t)\|_{L^p}<\infty,\qquad P-a.s.$

\noindent{\bf Theorem 2.1.3} Let $G$ satisfy the Lipschitz condition
$\| G(u)-G(v)\|_{L^2(K,H)}\leq \beta|u-v|$
for all $u,v\in dom(G)$, for some $\beta\in \mathbb{R}$ independent of $u,v$.
Then (2.1) admits at most one (probabilistically strong) solution such that $\sup_{t\in [0,T]}\|\Lambda^{1-\alpha+\varepsilon}\theta(t)\|_{L^p}<\infty, P-a.s.$ with $\frac{1}{p}\leq\frac{\alpha+\varepsilon}{2}$ and $\varepsilon\in (0,\alpha]$.

\noindent{\bf 2.2. The subcritical case}-----In this section, we will consider the subcritical case.

\noindent{\bf Theorem 2.2.1} Assume $\alpha>1/2$ and that $G$ does not depend on $\theta$ with $\rm{Tr}(\Lambda^{2-2\alpha+\varepsilon}GG^*)<\infty$ for some $\varepsilon>0$. Then for each initial condition $x\in H$, there exists a (probabilistically strong) solution $\theta$ to (2.1) over $[0,T]$ with initial condition $\theta(0)=x$.

\noindent{\bf Theorem 2.2.2} Assume $\alpha>\frac{1}{2}$. If  $G$ satisfies the following condition
$$\| \Lambda^{-1/2}(G(u)-G(v))\|_{L_2(K,H)}\leq \beta|\Lambda^{-1/2}(u-v)|\eqno(2.5)$$
for all $u,v\in dom(G)$, for some $\beta\in \mathbb{R}$ independent of $u,v$, then (2.1) admits at most one (probabilistically strong) solution $\theta$ such that $\sup_{t\in [0,T]}\|\theta(t)\|_{L^q}<\infty, P-a.s.$ for some $q$ with $0\leq 1/q<\alpha-\frac{1}{2}$

\noindent{\bf Corollary 2.2.3} Assume $\alpha>\frac{1}{2}$. If there exists a (probabilistically strong) solution  $\theta$ such that $\sup_{t\in [0,T]}\|\theta(t)\|_{L^q}<\infty,P-a.s.$ for some $q$ with $0\leq 1/q<\alpha-\frac{1}{2}$  and $G$ satisfies (2.5), then (2.1) admits only one such solution.

\noindent{\bf Theorem 2.2.4}  Assume $\alpha>\frac{1}{2}$ and $G$ satisfies (2.4), (2.5) and (G.1). Then for each initial condition $x\in H^1$, there exists a pathwise unique (probabilistic strong) solution $\theta$ of equation (2.1) over $[0,T]$ with  initial condition $\theta(0)=x$. Moreover, the solution satisfies $\sup_{t\in [0,T]}\|\theta(t)\|_{L^p}<\infty,P-a.s.$

\noindent{\bf Theorem 2.2.5}  Assume $\alpha>\frac{1}{2}$ and that  $G\in L^2(K,H)$ does not depend on $\theta$.
 Then for each initial condition $x\in L^p$ for some $p$ with $0\leq 1/p<\alpha-\frac{1}{2}$, there exists a pathwise unique (probabilistically strong)  solution $\theta$ of equation (2.1) over $[0,T]$ with initial condition $\theta(0)=x$. Moreover, this solution satisfies $\sup_{t\in [0,T]}\|\theta(t)\|_{L^p}<\infty,P-a.s.$

\noindent{\bf Theorem 2.2.6} (Markov property) Assume $\alpha>\frac{1}{2}$ and that $G\in L^2(K,H)$ does not depend on $\theta$. If $x\in L^p$ for some $p$ with $0\leq 1/p<\alpha-\frac{1}{2}$, for every bounded, $\mathcal{B}(H)$-measurable $F:H\rightarrow \mathbb{R}$, and all $s,t\in [0,T]$, $s\leq t, E(F(\theta(t))|\mathcal{F}_s)(\omega)=E(F(\theta(t,s,\theta(s)(\omega)))) \textrm{ for } P-a.e. \omega\in \Omega,$
where $\theta(t,s,\theta(s)(\omega))$ denotes the solution starting from $\theta(s)$ at time $s$.

Set $p_t(x,dy):=P\circ (\theta(t,x))^{-1}(dy), 0\leq t\leq T, x\in H,$
and for $\mathcal{B}(H)$-measurable $F:H\rightarrow \mathbb{R}$, and $t\in [0,T], x\in H, P_tF(x):=\int F(y)p_t(x,dy),$
provided $F$ is $p_t(x,dy)$-integrable. Then by Theorem 2.2.6, we have for $F:H\rightarrow \mathbb{R}$, bounded and $\mathcal{B}(H)$-measurable, $s,t\geq0, P_s(P_tF)(x)=P_{s+t} F(x), x\in L^p \textrm{ for some }p \textrm{ with }0\leq 1/p<\alpha-\frac{1}{2}.$

\noindent{\bf 3. Ergodicity and Exponential convergence for $\alpha>\frac{2}{3}$ }

\noindent{\bf Assumption 3.1} There are an isomophism $Q_0$ of $H$ and a number $s\geq1$ such that $G=A^{\frac{-s-\alpha}{2\alpha}}Q_0^{1/2}$.

 Set $\mathcal{W}=D(\Lambda^s)$ and  $|x|_\mathcal{W}=|\Lambda^s  x|$. Then by a similar method as in \cite{FR} and using the abstract results in \cite{GM} for exponential convergence, we obtain the following results.

\noindent{\bf Theorem 3.2} Assume $\alpha>\frac{2}{3}$ and Assumption 3.1. Then there exists a unique invariant measure $\nu$ on $\mathcal{W}$ for the transition semigroup $(P_t)_{t\geq0}$. Moreover:

(i) The invariant measure $\nu$ is ergodic in the sense of \cite{DZ96}.

(ii) The transition semigroup $(P_t)_{t\geq0}$ is $\mathcal{W}$-strong Feller, irreducible, and therefore strongly mixing in the sense of \cite{DZ96}.

(iii) There is $C_{\exp}>0$ and $a>0$ such that
$\|P_t^*\delta_{x_0}-\mu\|_{TV}\leq C_{\exp}(1+|x_0|^2)e^{-at},$
for all $t>0$ and $x_0\in H$, where $\|\cdot\|_{TV}$ is the total variation distance for measures.

\noindent{\bf 4. Markov Selections in the general case}-----By using the abstract results for Markov selections in \cite{GRZ}, we obtain the following results.

\noindent{\bf Theorem 4.1} Assume $G$ satisfies (G.1). Then there exists an almost sure Markov family $(P_{x_0})_{x_0\in H}$ for Eq. (2.1).

\noindent{\bf Acknowledgement} We thank Wei Liu for very helpful discussions.

\end{document}